\documentclass[mathpazo]{cicp}

\usepackage{graphicx}
\usepackage{bm}
\usepackage{color}
\usepackage{booktabs}
\usepackage{multirow}

\begin{document}
\title{A filtered MgNet solver for radiative transfer equations}

\author{Qinchen Song\affil{1}, Lei Zhang\affil{1}\corrauth, Xinliang Liu\affil{2}\corrauth}
\address{\affilnum{1}\ Department of Mathematics, Institute of Natural Sciences,
Shanghai Jiao Tong University, Minhang, Shanghai 200240, P.R. China.\\
\affilnum{2}\ Extreme Computing Research Center, King Abdullah University of Science and Technology}
\email{{\tt sqc9931@sjtu.edu.cn}, {\tt lzhang2012@sjtu.edu.cn}, {\tt xinliang.liu@kaust.edu.sa}}


\begin{abstract}
Conventional numerical solvers for the radiative transfer equation (RTE) exhibit severe sensitivity to medium parameters. To address this, we propose an operator learning framework that approximates the RTE solution map as a function of material properties. The core architecture, MgNet, preserves the solution operator framework established by recursive skeleton factorization (RSF) but substitutes its coefficient-specific sub-operators (e.g. smoother, prolongation operator and restriction operator) with learnable neural components. This design transcends the the fixed parametric structure of classical schemes, enabling data-driven sub-operator optimization and learning of their medium-parameter dependence. To mitigate spectral bias in operator learning, we introduce an adaptive angular compression technique within the loss function that dynamically suppresses high-frequency modes responsible for training instability. Comprehensive benchmarks demonstrate that, when deployed as a learned preconditioner, MgNet achieves at least 10 times acceleration over conventional preconditioners in the diffusive regime and maintains robust generalization to unseen parameter configurations. By unifying multilevel factorization structure with deep operator learning, this work establishes a physics-constrained operator-learning paradigm for radiative transport simulations.
\end{abstract}

\ams{35Q70, 65M55, 65M50, 65M06}
\keywords{radiative transport equation (RTE); discrete ordinates method (DOM); adaptive tailored finite point scheme (ATFPS); spectral bias mitigation; multilevel method; MgNet; operator learning.}

\maketitle

\section{Introduction}
The radiative transfer equation (RTE) governs the propagation and scattering of particles (e.g., photons, neutrons) in heterogeneous media, with foundational applications in fields such as nuclear engineering \cite{wang2024simulation,connolly2013heterogeneous}, atmospheric physics \cite{zhang2024advances}, astrophysics \cite{chandrasekhar2013radiative,macdonald2022trident}, remote sensing \cite{mishchenko2003radiative} and optical tomography \cite{klose2002optical}. The steady state RTE under diffusive scaling reads
\begin{equation}\label{eq:RTE}
    \mathbf{u}\cdot\nabla \psi(\mathbf{z},\mathbf{u})+\frac{\sigma_{T}}{\epsilon}(\mathbf{z})\psi(\mathbf{z},\mathbf{u})=\Big(\frac{\sigma_{T}(\mathbf{z})}{\epsilon}-\epsilon\sigma_{a}(\mathbf{z})\Big)\int_{\mathbb{S}^{2}}\kappa(\mathbf{u},\mathbf{u}')\psi(\mathbf{z},\mathbf{u})\,d\mathbf{u}'+\epsilon q(\mathbf{z}), 
\end{equation}
where $\mathbf{z}\in \Omega \subset\mathbb{R}^{3}$ and  $\mathbf{u}\in\mathbb{S}^{2}$ denote spatial position, and particle direction, respectively. The angular flux $\psi(\mathbf{z},\mathbf{u})$ represents the particle density moving along $\mathbf{u}$ at $(\mathbf{z})$. The dimensionless parameter $\epsilon \in (0,1]$ scales the mean free path relative to the domain characteristic length, governing the transition between the transport regime ($\epsilon \sim \mathcal{O}(1)$) and the diffusive regime ($\epsilon \ll 1$). Following the diffusive scaling framework in \cite{han2014two}, the total cross section, scattering cross section, absorption cross section, and external source are rescaled to preserve asymptotic consistency as $\epsilon \to 0$. We assume $\sigma_T$, $\sigma_a$, and $q$ are non-negative, spatially dependent, piecewise smooth, and bounded, with uniformly bounded gradients except at material interfaces. This regularity ensures that standard spatial discretizations adequately resolve material variations.

The scattering kernel $\kappa(\mathbf{u},\mathbf{u}')$ models the probability of directional transition upon collision. A widely used choice is the Henyey--Greenstein (HG) kernel \cite{henyey1941diffuse},
\begin{equation}\label{eq:hg}
\kappa(\mathbf{u},\mathbf{u}')=\frac{1-g^{2}}{(1+g^{2}-2g\,\mathbf{u}\cdot\mathbf{u}')^{3/2}},
\end{equation}
where $g\in[-1,1]$ is the anisotropy factor. The system is closed by the 
inflow Dirichlet boundary conditions
\begin{equation}\label{eq:boundary_condition}
    \psi(\mathbf{z},\mathbf{u})=\Psi(\mathbf{z},\mathbf{u}),\quad (\mathbf{z},\mathbf{u})\in\Gamma^{-}=\{(\mathbf{z},\mathbf{u})\in\partial\Omega\times\mathbb{S}^{2}\mid \mathbf{u}\cdot\mathbf{n}_{\mathbf{z}}<0 \},  
\end{equation}
and interface continuity
\begin{equation}\label{eq:interface_condition}
    [\psi]_{\alpha}=0,
\end{equation}
where $[\cdot]_{\alpha}$ denotes the jump across interface $\alpha$ and $\mathbf{n}_{\mathbf{z}}$ is the outward unit normal.

Discretizing the RTE yields high-dimensional, strongly coupled linear systems whose size scales with spatial, angular, and temporal resolution. The multiscale nature induced by $\epsilon$ exacerbates this challenge: in the diffusive regime, boundary layers and sharp gradients demand highly refined meshes to maintain asymptotic consistency and numerical stability \cite{lewis1984computational,larsen1987asymptotic,adams2001discontinuous}. Consequently, classical solvers face a severe computational bottleneck. Direct factorization methods (e.g., LU, Cholesky) offer robustness but incur $O(N^3)$ complexity and $O(N^2)$ memory, rendering them impractical for large-scale problems. While hierarchical low-rank techniques such as recursive skeleton factorization (RSF) mitigate this scaling \cite{ho2012fast}, their factorization remains tightly coupled to medium coefficients. Iterative Krylov methods (e.g., GMRES, BiCGSTAB) scale more favorably but rely heavily on preconditioners whose effectiveness deteriorates rapidly under parameter variations or regime transitions \cite{saad2003iterative}. In heterogeneous or time-varying settings, this strong parameter dependence necessitates frequent recomputation of factorizations or preconditioner updates, incurring substantial overhead and limiting scalability.

To overcome this rigidity, we apply a parameter-aware operator learning framework that directly approximates the RTE solution map as a function of material properties. The core architecture, MgNet \cite{he2019mgnet,zhu2024enhanced}, preserves the hierarchical solution-operator topology of classical RSF structure but replaces its fixed, coefficient-dependent sub-operators with learnable neural components. This design transcends the rigid parametrization of classical schemes, enabling data-driven optimization of transport operators while explicitly encoding their nonlinear dependence on medium parameters. The trained network serves as a robust, parameter-invariant preconditioner within iterative solvers.

Besides, training is driven by an unsupervised physics‑informed loss constructed from the discrete linear‑system residual \cite{raissi2019physics,cuomo2022scientific,plankovskyy2025review}, eliminating the need for labeled solutions and reducing computational cost while embedding physical consistency to improve generalization. To overcome the spectral bias that hinders high‑frequency angular modes, we incorporate an adaptive angular compression strategy \cite{song2025adaptive} into the loss, which dynamically suppresses stiff components to stabilize convergence and accelerate training.


The application of deep learning to the radiative transfer equation (RTE) has evolved from simple data-driven emulation to physics-constrained operator learning. Early supervised approaches, such as those by Lagerquist et al. \cite{lagerquist2021using} and Brodrick et al. \cite{brodrick2021generalized}, utilized CNNs and feed-forward networks to accelerate retrievals. While efficient, these methods are constrained by the need for massive labeled datasets. To mitigate this, Zhu et al. \cite{zhu2026deeprte} proposed DeepRTE using an attention-based pre-training strategy.

Alternatively, Physics-Informed Neural Networks (PINNs) \cite{mishra2021physics} offer an unsupervised framework by embedding the RTE residual into the loss function. This paradigm has been extended to handle heterogeneous media via Fourier features \cite{huhn2023physics}, boundary dependencies \cite{xie2024boundary}, and multiscale asymptotic limits \cite{jin2023asymptotic, lu2022solving}. Despite their success, standard PINNs often struggle with the "spectral bias," making it difficult to resolve high-frequency layer information characteristic of transport solutions.

To enhance multiscale feature extraction, the MgNet framework \cite{he2019mgnet, zhu2024enhanced} unifies multigrid methods with iterative network structures. However, applying MgNet within a physics-informed context—specifically through specialized spatial discretizations designed to actively filter high-frequency layer components—remains an open challenge that this work aims to address.

The principal contributions of this work are:
\begin{itemize}
    \item \textbf{Interpretable MgNet Design:} We design an MgNet architecture whose structure and hyperparameters (e.g., channel width) are theoretically aligned with the RTE solution operator's structure \cite{song2026fast}.
    \item \textbf{Layer-Filtering Loss:} By embedding a compressed spatial discretization \cite{song2025adaptive} into the physics-informed loss, we explicitly remove high-frequency modes associated with layers, significantly enhancing training convergence and robustness.
\end{itemize}

The remainder of this paper is organized as follows. Section \ref{sec:discretization_scheme} reviews the discrete formulation of the RTE and introduces the adaptive angular domain compression technique, which is embedded into the spatial discretization. Section \ref{sec:MgNet} details the MgNet architecture and the physics-informed loss construction. Section \ref{sec:numerical_experiment} presents numerical benchmarks across diffusive and transport regimes, reporting training/validation losses, iterative performance as preconditioner, with comparisons to classical preconditioner on random test cases. Finally, Section \ref{sec:summary} concludes the paper and discusses possible extensions.

\section{Discretization scheme}
\label{sec:discretization_scheme}
In the following sections, we focus on the solution of the steady-state RTE in x-y geometry, where the spatial domain is reduced to two dimensions under the assumption of invariance in the z-direction. Accordingly, the original RTE \eqref{eq:RTE} reduces to:
\begin{equation}\label{eq:RTE_xy}
    \begin{aligned}
    &c\partial_{x}\psi(x,y,c,s)+s\partial_{y}\psi(x,y,c,s) +\frac{\sigma_{T}(x,y)}{\epsilon}\psi(x,y,c,s)\\
    &=\Big(\frac{\sigma_{T}(x,y)}{\epsilon}-\epsilon\sigma_{a}(x,y)\Big)\int_{\mathbb{D}^{2}}\kappa(c,s,c',s')\psi(x,y,c',s')dc'ds'+\epsilon q(x,y), \\
    \end{aligned}
\end{equation}
where $c$ and $s$ denote the projections of the unit direction vector $\mathbf{u}$ onto the $x$- and $y$-axes, respectively, and $\mathbb{D}^{2} = \{(x,y) \mid x^2 + y^2 \leq 1\}$ represents the projection of the unit sphere $\mathbb{S}^2$ onto the $x$-$y$ plane.

In this section, we present the discretization of the steady-state RTE in x-y geometry using DOM \cite{lewis1984computational} for the angular domain and TFPS \cite{han2014two} for the spatial domain. These discretizations enable the construction of a physics-informed loss function that enforces compliance with the discretized RTE. Additionally, we introduce an adaptive angular domain compression technique, which combined with TFPS yields the adaptive tailored finite point scheme (ATFPS)\cite{song2025adaptive}. Embedding this compression into the loss function filters high-frequency modes, thus accelerating convergence and enhancing robustness.

\subsection{discrete ordinate method (DOM)}
\label{subsec:DOM}
\begin{figure}[htbp]
    \centering
    \includegraphics[width=0.7\textwidth]{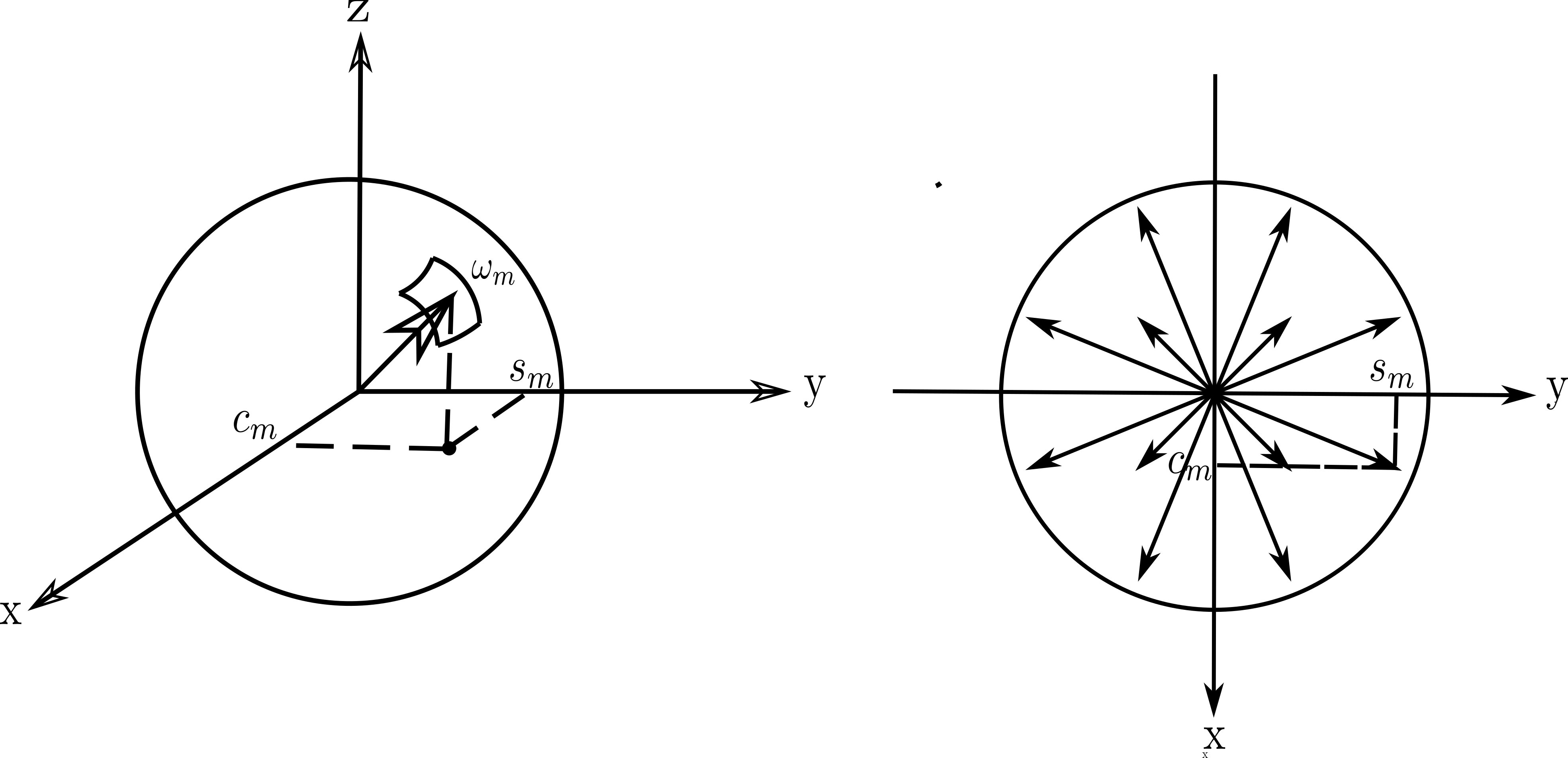}
    \caption{DOM in x-y geometry. Left figure: a quadrature point in three spatial dimensional case. Right figure: example of DOM ($S_{N}$) in x-y geometry with $N=2$.}
    \label{fig:DOM}
\end{figure}

The DOM quadrature points in x-y geometry are obtained by projecting their three-dimensional counterparts onto the two-dimensional plane, as illustrated in Figure \ref{fig:DOM}. We denote the quadrature set for DOM in x-y geometry as $\{\mathbf{u}_{m},\omega_{m}\}_{m\in\mathcal{M}}$ with $\mathcal{M}=\{1,2,\dots,4M\}$ being the index set for discrete velocity directions. Here $\mathbf{u}_{m}=(c_{m},s_{m})=((1-\zeta_{m}^{2})^{1/2}\cos(\theta_{m}), (1-\zeta_{m}^{2})^{1/2}\sin(\theta_{m}))$ represents the $m$-th velocity direction with $\zeta_{m}\in[0,1]$ and $\theta_m\in[0,2\pi)$, and $\omega_{m}$ represents the corresponding weight. 

We approximate the integral term in the RTE by its numerical quadrature, and denote the approximation of $\psi(x,y,c_{m},s_{m})$ as $\psi_{m}(x,y)$. Then we get the following 2D  \textit{discrete ordinate RTE} in x-y geometry:
\begin{equation}\label{eq:DORTE_2d}
    c_{m}\frac{\partial}{\partial x}\psi_{m}(x,y)+s_{m}\frac{\partial}{\partial y}\psi_{m}(x,y)+\sigma_{T}\psi_{m}(x,y)=\sigma_{s}\sum_{n\in V}\kappa_{mn}\psi_{n}(x,y)\omega_{n}+ q,
\end{equation}
for $m\in\mathcal{M}$. Here $\kappa_{mn}$ is an approximation of $\kappa(c_{m},s_{m},c_{n}, s_{n})$ \cite{chen2018uniformly}. In the following, we let $\bm{\psi}(x,y)$ be the $4M$-dimensional vector function $\Big(\psi_{1}(x,y),\psi_{2}(x,y),\dots,\psi_{4M}(x,y)\Big)^{T}$.

For simplicity, we consider the spatial domain $\Omega = [0,1] \times [0,1]$. 
The boundary conditions are discretized as:
\begin{equation}
\label{eq:boundary_condition_2d}
\begin{aligned}
\psi_{m}(0,y) = \Psi_{m}(0,y),\quad \forall y\in[0,1],\ c_{m}>0; & \quad \psi_{m}(1,y) = \Psi_{m}(1,y),\quad \forall y\in[0,1],\ c_{m}<0; \\
\psi_{m}(x,0) = \Psi_{m}(x,0),\quad \forall x\in[0,1],\ s_{m}>0; & \quad \psi_{m}(x,1) = \Psi_{m}(x,1),\quad \forall x\in[0,1],\  s_{m}<0.
\end{aligned}
\end{equation}
At material interfaces $\alpha$, the flux continuity condition reads:
\begin{equation}
\label{eq:interface_condition_2d}
[\psi]|_{\alpha} = 0.
\end{equation}

\subsection{tailored finite point scheme (TFPS)}
\label{subsec:TFPS}
Assuming a uniform mesh and that parameter discontinuities occur only at cell interfaces, we denote grid points $x_i = ih, y_j = jh$ ($0 \leq i,j \leq I$) with $h = 1/I$. The cell $C_{i,j} = [x_{i-1}, x_i] \times [y_{j-1}, y_j]$ and the collection of all cells is $\mathfrak{C} = \{C_{i,j}\}_{1 \leq i,j \leq I}$. For any cell $C \in \mathfrak{C}$, we write $C = [x_{C,l}, x_{C,r}] \times [y_{C,b}, y_{C,t}]$. The cell averages are defined as
\[
\sigma_{T,C}=\frac{1}{|C|}\int_C\sigma_{T}(x,y)\,dx\,dy,\quad \sigma_{s,C}=\frac{1}{|C|}\int_C\sigma_{s}(x,y)\,dx\,dy,\quad q_{C}=\frac{1}{|C|}\int_C q(x,y)\,dx\,dy,
\]
where $|C|$ denotes the cell area. We construct piecewise constant functions $\bar{\sigma}_{T}, \bar{\sigma}_{s}, \bar{q}$ that take values $\sigma_{T,C}, \sigma_{s,C}, q_{C}$ on $C$, respectively. For simplicity, we assume $\sigma_{a,C} = \sigma_{T,C} - \sigma_{s,C} \neq 0$ for all cells $C$.

Using these piecewise constant approximations, Eq. \eqref{eq:DORTE_2d} becomes
\begin{equation}
\label{discretization_2D_TFPS}
c_{m}\partial_{x}\psi_{m}+s_{m}\partial_{y}\psi_{m}+\bar{\sigma}_{T}\psi_{m}=\bar{\sigma}_s\sum_{n\in\mathcal{M}}\kappa_{m,n}\psi_{n}\omega_{n}+\bar{q},\quad m\in\mathcal{M}.
\end{equation}
Within any cell $C \in \mathfrak{C}$, Eq. \eqref{discretization_2D_TFPS} reduces to a constant-coefficient equation:
\begin{equation}
\label{discretization_2D_TFPS_local}
c_m \partial_{x} \psi_m + s_m \partial_{y} \psi_m + \sigma_{T,C} \psi_m = \sigma_{s,C} \sum_{n \in \mathcal{M}} \kappa_{m,n} \psi_n \omega_n + q_C,\quad m \in \mathcal{M}.
\end{equation}

Consider the corresponding homogeneous equation:
\begin{equation}
\label{discretization_2D_TFPS_local_homo}
c_m \partial_{x} \psi_m + s_m \partial_{y} \psi_m + \sigma_{T,C} \psi_m = \sigma_{s,C} \sum_{n \in \mathcal{M}} \kappa_{m,n} \psi_n \omega_n,\quad m \in \mathcal{M}.
\end{equation}
Clearly, the solution of \eqref{discretization_2D_TFPS_local_homo} admits the following form:
\begin{equation}
\label{eq:2D_eigenfunc_pre}
\psi_{m}(x,y) = l_{m}\exp\left\{\sigma_{T,C}(\xi_{x}x+\xi_{y}y)\right\}.
\end{equation}
Substituting \eqref{eq:2D_eigenfunc_pre} into \eqref{discretization_2D_TFPS_local_homo} yields the eigenvalue problem:
\begin{equation}
\label{discretization_2D_eigenfunc_constraint}
\left[\frac{\sigma_{s,C}}{\sigma_{T,C}}KW-I\right]l = (\xi_{x} C + \xi_{y}S)l,
\end{equation}
where $C = \operatorname{diag}\{c_m\}_{m\in\mathcal{M}}$, $S = \operatorname{diag}\{s_m\}_{m\in\mathcal{M}}$, $W = \operatorname{diag}\{\omega_m\}_{m\in\mathcal{M}}$, and $K = (\kappa_{m,n})_{m,n\in\mathcal{M}}$.

There are infinitely many eigenpairs $(\xi_{x},\xi_{y},l)$ satisfying \eqref{discretization_2D_eigenfunc_constraint}. To construct basis functions, we consider two special cases:
\begin{enumerate}
\item Setting $\xi_{y}=0$, the eigenpairs $(\xi_{x},l)$ correspond to the matrix $C^{-1}[\frac{\sigma_{s,C}}{\sigma_{T,C}}KW-I]$. There are $4M$ eigenpairs, denoted $(\xi_{x,C}^{k},l_{x,C}^{k})_{k=-2M,\dots,-1,1,\dots,2M}$, ordered as
\[
\xi_{x,C}^{-2M}<\cdots<\xi_{x,C}^{-1}< 0 < \xi_{x,C}^{1}<\cdots<\xi_{x,C}^{2M}.
\]
\item Setting $\xi_{x}=0$, the eigenpairs $(\xi_{y},l)$ correspond to the matrix $S^{-1}[\frac{\sigma_{s,C}}{\sigma_{T,C}}KW-I]$, yielding $4M$ eigenpairs $(\xi_{y,C}^{k},l_{C}^{k})_{k=-2M,\dots,-1,1,\dots,2M}$ ordered as
\[
\xi_{y,C}^{-2M}<\cdots<\xi_{y,C}^{-1}< 0< \xi_{y,C}^{1}<\cdots<\xi_{y,C}^{2M}.
\]
\end{enumerate}
All eigenvectors are normalized such that $\| l_{x,C}^{k} \|_{\infty} =\| l_{y,C}^{k} \|_{\infty} = 1$.

The corresponding basis functions are:
\begin{equation}
\label{eq:2D_eigenfunc}
\chi_{C}^{k}(\mathbf{z}) = l_{C}^{k}\exp\left\{\sigma_{T,C}(\bm{\xi}_{C}^{k})^{T}(\mathbf{z}-\mathbf{z}_{C}^{k})\right\},\quad \mathbf{z}\in C,\quad 1\leq \vert k\vert \leq 4M,
\end{equation}
where $\mathbf{z}=(x,y)$. The eigenvalue vector $\bm{\xi}_{C}^{k}$, eigenvector $l_{C}^{k}$, and reference point $\mathbf{z}_{C}^{k}$ are defined as:
\[
\bm{\xi}_{C}^{k} =
\begin{cases}
(\xi_{x,C}^{k+2M},0)^{T}, & -4M\leq k\leq -1,\\
(0,\xi_{y,C}^{k-2M})^{T}, & 1\leq k\leq 4M,
\end{cases}
\quad
l_{C}^{k} =
\begin{cases}
l_{x,C}^{k+2M}, & -4M\leq k\leq -1,\\
l_{y,C}^{k-2M}, & 1\leq k\leq 4M,
\end{cases}
\]
\[
\mathbf{z}_{C}^{k} =
\begin{cases}
(x_{C,l}, y_{C,c}), & -4M\leq k\leq -2M-1,\\
(x_{C,r}, y_{C,c}), & -2M\leq k\leq -1,\\
(x_{C,c}, y_{C,b}), & 1\leq k\leq 2M,\\
(x_{C,c}, y_{C,t}), & 2M+1\leq k\leq 4M,
\end{cases}
\]
where $x_{C,c}$ and $y_{C,c}$ are the $x$- and $y$-coordinates of the cell center. The choice of reference points ensures that eigenfunctions remain bounded by 1 within each cell, enhancing numerical stability.

The general solution of the homogeneous equation \eqref{discretization_2D_TFPS_local_homo} spanned by these basis functions is:
\begin{equation}
\label{discretization_2D_local_sol_homo1}
\bm{\psi}(\mathbf{z}) = \sum_{k\in\mathcal{V}}\alpha_{C}^{k}\chi_{C}^{k}(\mathbf{z}),\quad \mathcal{V} = \{-4M,\dots,-1, 1,\dots,4M\}.
\end{equation}
A particular solution of the inhomogeneous equation \eqref{discretization_2D_TFPS_local} is $\chi_{C}^{s} = \frac{q_C}{\sigma_{T,C} - \sigma_{s,C}} \mathbf{1}$. Consequently, the global solution of \eqref{discretization_2D_TFPS} can be expressed as:
\begin{equation}
\label{eq:2DTFPS_sol}
\bm{\psi}(\mathbf{z}) = \sum_{C\in\mathfrak{C}}\mathbb{I}_C(\mathbf{z})\left(\sum_{k\in\mathcal{V}}\alpha_{C}^{k}\chi_{C}^{k}(\mathbf{z}) + \chi_{C}^{s}\right),
\end{equation}
where $\mathbb{I}_C$ is the indicator function of cell $C$.

Let $\mathfrak{I}$ denote the set of all cell interfaces, with $\mathfrak{I}_{\mathrm{b}}$ and $\mathfrak{I}_{\mathrm{in}}$ representing the physical boundary and interior interfaces, respectively. For an interior interface $\mathfrak{i}=C_-\cap C_+\in \mathfrak{I}_{\mathrm{in}}$, the continuity condition $\bm{\psi}|_{C_{-}}(\mathbf{z})=\bm{\psi}|_{C_{+}}(\mathbf{z})$ for all $\mathbf{z}\in\mathfrak{i}$ should theoretically hold. To formulate a well-posed linear system, TFPS relaxes this to continuity at the midpoint $\mathfrak{i}_{\mathrm{mid}}$:
\begin{equation}
\label{eq:2DTFPS_continuity_cond}
\sum_{k\in\mathcal{V}}\alpha_{C_{-}}^{k}\chi_{C_{-}}^{k}(\mathfrak{i}_{\mathrm{mid}}) + \chi_{C_{-}}^{s} = \sum_{k\in\mathcal{V}}\alpha_{C_{+}}^{k}\chi_{C_{+}}^{k}(\mathfrak{i}_{\mathrm{mid}}) + \chi_{C_{+}}^{s},\quad \forall \mathfrak{i}\in\mathfrak{I}_{\mathrm{in}}.
\end{equation}
For a boundary interface $\mathfrak{i}=C\cap\partial\Omega\in\mathfrak{I}_{\mathrm{b}}$, the incoming flux at the midpoint must satisfy the boundary condition:
\begin{equation}
\label{eq:2DTFPS_boundary_cond}
\sum_{k\in\mathcal{V}}\alpha_{C}^{k}\chi_{C,\mathfrak{i}-}^{k}(\mathfrak{i}_{\mathrm{mid}}) + \chi_{C,\mathfrak{i}-}^{s} = \bm{\psi}_{\mathfrak{i}-}(\mathfrak{i}_{\mathrm{mid}}).
\end{equation}
Bold symbols here denote vectors over all angular directions, e.g., 
$\bm{\psi}(x,y)=\big(\psi_{m}(x,y)\big)_{m\in\mathcal{M}_{\mathrm{2d}}}$. 
The subscript "$\mathfrak{i}-$" indicates restriction to incoming directions 
at interface $\mathfrak{i}$, i.e., 
$\chi_{C,\mathfrak{i}-}^{k} = (\chi_{C,m}^{k})_{\mathbf{\Omega}_{m}\cdot \mathbf{n}_{\mathfrak{i}}<0}$ 
and $\bm{\psi}_{\mathfrak{i}-} = (\psi_{m})_{\mathbf{\Omega}_{m}\cdot\mathbf{n}_{\mathfrak{i}}<0}$, 
where $\mathbf{n}_{\mathfrak{i}}$ is the outward normal.

Taking the basis function coefficients $\alpha = \{\alpha_{C}^{k}\}_{C\in\mathfrak{C}, k\in\mathcal{K}}$ as unknowns, the linear system is constructed from the continuity conditions \eqref{eq:2DTFPS_continuity_cond} and boundary conditions \eqref{eq:2DTFPS_boundary_cond}:
\begin{equation}
\label{eq:2DTFPS_linear_sys}
A\alpha = b,
\end{equation}
where $\alpha = (\alpha_{C}^{k})_{k\in\mathcal{V},C\in\mathfrak{C}_{\mathrm{2d}}}$, matrix $A$ consists of basis function values at interface midpoints, and the right-hand side $b$ contains boundary values and particular solutions. Once the linear system \eqref{eq:2DTFPS_linear_sys} is solved for $\alpha$, substituting back into \eqref{eq:2DTFPS_sol} yields the numerical solution of the two-dimensional discrete ordinates radiative transfer equation \eqref{eq:DORTE_2d}.

\subsection{adaptive tailored finite point scheme}
\label{subsec:ATFPS}
As demonstrated in the asymptotic analysis of TFPS \cite{chen2018uniformly}, in special cases such as the diffusion regime or sharp interface problems where transport and diffusion coexist, a reduced set of TFPS basis functions (specifically the slowly decaying ones) can be employed to derive a compressed linear system.

Inspired by this observation, we note that the rapidly decaying basis functions primarily capture information associated with strongly varying structures such as boundary and interface layers, while the slowly decaying basis functions correspond to the smooth components of the solution. These two components can be approximately decoupled. Accordingly, we employ only the slowly decaying basis functions to derive the smooth part of the solution and subsequently reconstruct the layer information to get a solution with global accuracy. This constitutes the adaptive angular domain compression technique. Integrating this technique with TFPS yields the adaptive tailored finite point scheme (ATFPS).

The remainder of this subsection is organized as follows. We first introduce the ATFPS basis, which corresponds to the slowly decaying portion of the TFPS basis, along with necessary notation. We then derive the corresponding boundary and interface conditions and present the compressed linear system. Finally, we describe the reconstruction technique for recovering layer information.

\subsubsection{Basis functions and solution representation}
\label{subsubsec:ATFPS_basis}

Specifically, we first introduce a tolerance threshold $\delta$. Based on $\delta$ and the optical properties of the cell (such as the scattering ratio $\gamma = \sigma_s / \sigma_T$, the anisotropy factor $g$, and the optical thickness $\sigma_T h$), the ATFPS selects a subset of the TFPS basis $\{\chi_{C}^{k}\}_{k\in\mathcal{V},C\in\mathfrak{C}}$ to form a new set of smoother basis functions $\{\chi_{C}^{k}\}_{k\in\mathcal{V}_{\delta,C},C\in\mathfrak{C}}$. The index set $\mathcal{V}_{\delta,C}$ is defined as:
\begin{equation}
\label{eq:def_V_delta}
\mathcal{V}_{\delta,C} = \big\{ k \mid \|\chi_{C}^{k}(\mathbf{x}_C)\|_{\infty} > \delta \big\}
= \big\{ k \mid \exp\{-\tfrac{1}{2}|\lambda_{C}^{k}|\sigma_{T,C}h\} > \delta \big\},
\end{equation}
where $\mathbf{x}_C$ denotes the cell center.

Furthermore, we define the complement $\bar{\mathcal{V}}_{\delta,C} = \mathcal{V} \setminus \mathcal{V}_{\delta,C}$. Thus, for any cell $C \in \mathfrak{C}$, the basis functions $\{\chi_{C}^{k}\}_{k\in\mathcal{V}}$ can be divided into two groups: $\{\chi_{C}^{k}\}_{k\in\mathcal{V}_{\delta,C}}$ (the selected adaptive TPFS basis functions with $\|\chi_{C}^{k}(\mathbf{x}_C)\|_{\infty} > \delta$) and $\{\chi_{C}^{k}\}_{k\in\bar{\mathcal{V}}_{\delta,C}}$ (the unselected rapidly decaying basis functions with $\|\chi_{C}^{k}(\mathbf{x}_C)\|_{\infty} < \delta$).

The smooth component $\bar{\bm{\psi}}_{\delta}$ of the adaptive TPFS solution can be spanned by the selected basis functions:
\begin{equation}
\label{eq:2DAdaptiveTFPS_sol_smooth}
\bar{\bm{\psi}}_{\delta} = \sum_{C\in\mathfrak{C}}\sum_{k\in\mathcal{V}_{\delta,C}}\alpha_{\delta,C}^{k}\chi_{C}^{k} + \chi_{C}^{s},
\end{equation}
while the layer component $\tilde{\bm{\psi}}_{\delta}$ is spanned by the unselected basis functions:
\begin{equation}
\label{eq:2DAdaptiveTFPS_sol_layer}
\tilde{\bm{\psi}}_{\delta} = \sum_{C\in\mathfrak{C}}\sum_{k\in\bar{\mathcal{V}}_{\delta,C}}\alpha_{\delta,C}^{k}\chi_{C}^{k}.
\end{equation}
Therefore, the overall solution of the adaptive TPFS can be expressed as:
\begin{equation}
\label{eq:2DAdaptiveTFPS_sol}
\bm{\psi}_{\delta} = \sum_{C\in\mathfrak{C}}\sum_{k\in\mathcal{V}}\alpha_{\delta,C}^{k}\chi_{C}^{k} + \chi_{C}^{s},
\end{equation}
where all coefficients $\alpha_{\delta,C}^{k}$ are to be determined. These coefficients can be determined from the subsequent continuity and boundary conditions.

\subsubsection{Newly defined notations}
\label{subsubsec:new_notation}

To facilitate the treatment of continuity and boundary conditions, we introduce some notation. For any cell $C \in \mathfrak{C}$, basis function index $k \in \mathcal{V}$, and interface $\mathfrak{i} \in C$, if the basis function $\chi_C^{k}$ attains its maximum on $\mathfrak{i}$, we say that this basis function is centered at $\mathfrak{i}$. Based on this, the index set $\mathcal{V}_C$ of local basis functions within cell $C$ can be divided into three parts:
\begin{itemize}
    \item $\mathcal{V}_{C,\mathfrak{i}}$: the index set of basis functions in $C$ centered at $\mathfrak{i}$;
    \item $\mathcal{V}_{C,\mathfrak{i}^{\parallel}}$: the index set of basis functions in $C$ centered at the interface opposite to $\mathfrak{i}$;
    \item $\mathcal{V}_{C,\mathfrak{i}^{\perp}}$: the index set of basis functions in $C$ centered at the interface perpendicular to $\mathfrak{i}$.
\end{itemize}

Based on the definitions of $\mathcal{V}_{\delta,C}$ and $\mathcal{V}_{C,\mathfrak{i}}$, the index set $\mathcal{V}_{\delta,C}$ of ATPFS basis functions within cell $C$ can be further partitioned into three parts centered at $\mathfrak{i}$, opposite to $\mathfrak{i}$, and perpendicular to $\mathfrak{i}$, denoted respectively as:
\[
\mathcal{V}_{\delta,C,\mathfrak{i}} = \mathcal{V}_{C,\mathfrak{i}} \cap \mathcal{V}_{\delta,C}, \quad
\mathcal{V}_{\delta,C,\mathfrak{i}^{\parallel}} = \mathcal{V}_{C,\mathfrak{i}^{\parallel}} \cap \mathcal{V}_{\delta,C}, \quad
\mathcal{V}_{\delta,C,\mathfrak{i}^{\perp}} = \mathcal{V}_{C,\mathfrak{i}^{\perp}} \cap \mathcal{V}_{\delta,C}.
\]

Similarly, the index set of unselected basis functions in cell $C$ can also be divided into the three categories, corresponding to being centered at $\mathfrak{i}$, opposite to $\mathfrak{i}$, and perpendicular to $\mathfrak{i}$, denoted as:
\[
\bar{\mathcal{V}}_{\delta,C,\mathfrak{i}} = \mathcal{V}_{C,\mathfrak{i}} \setminus \mathcal{V}_{\delta,C,\mathfrak{i}}, \quad
\bar{\mathcal{V}}_{\delta,C,\mathfrak{i}^{\parallel}} = \mathcal{V}_{C,\mathfrak{i}^{\parallel}} \setminus \mathcal{V}_{\delta,C,\mathfrak{i}^{\parallel}}, \quad
\bar{\mathcal{V}}_{\delta,C,\mathfrak{i}^{\perp}} = \mathcal{V}_{C,\mathfrak{i}^{\perp}} \setminus \mathcal{V}_{\delta,C,\mathfrak{i}^{\perp}}.
\]

Additionally, we define the index sets $\mathcal{V}_{\mathfrak{i}}$, $\mathcal{V}_{\delta,\mathfrak{i}}$, and $\bar{\mathcal{V}}_{\delta,\mathfrak{i}}$ as follows:
\[
\mathcal{V}_{\mathfrak{i}} = \bigcup_{\mathfrak{i} \cap C \neq \emptyset} \mathcal{V}_{C,\mathfrak{i}}, \quad
\mathcal{V}_{\delta,\mathfrak{i}} = \bigcup_{\mathfrak{i} \cap C \neq \emptyset} \mathcal{V}_{\delta,C,\mathfrak{i}}, \quad
\bar{\mathcal{V}}_{\delta,\mathfrak{i}} = \mathcal{V}_{\mathfrak{i}} \setminus \mathcal{V}_{\delta,\mathfrak{i}}.
\]
These represent, respectively, the index sets of TPFS basis functions, ATFPS basis functions, and their difference, all centered at interface $\mathfrak{i}$ within the physical domain.

For any interior interface $\mathfrak{i}=C_{-}\cap C_{+}\in\mathfrak{I}_{\mathrm{in}}$, 
let $\{\eta_{\mathfrak{i}}^{k}\}_{k\in\mathcal{V}_{\delta,\mathfrak{i}}}$ be the orthonormal 
basis obtained via QR decomposition of the set 
$\{l_{C_{-}}^{(k_{-})},l_{C_{+}}^{(k_{+})}\}_{k_{-}\in\mathcal{V}_{\delta,C_{-},\mathfrak{i}},k_{+}\in\mathcal{V}_{\delta,C_{+},\mathfrak{i}}}$. 
Define $\{\eta_{\mathfrak{i}}^{k}\}_{k\in\bar{\mathcal{V}}_{\delta,\mathfrak{i}}}$ as a reordering of 
$\{l_{C_{-}}^{(k_{-})},l_{C_{+}}^{(k_{+})}\}_{k_{-}\in\bar{\mathcal{V}}_{\delta,C_{-},\mathfrak{i}},k_{+}\in\bar{\mathcal{V}}_{\delta,C_{+},\mathfrak{i}}}$.  
For a boundary interface $\mathfrak{i}=C\cap\partial\Omega\in\mathfrak{I}_{\mathrm{b}}$, 
to match incoming boundary conditions, we take 
$\{\eta_{\mathfrak{i}}^{k}\}_{k\in\mathcal{V}_{\delta,\mathfrak{i}}}$ as the orthonormal basis 
from QR decomposition of $\{l_{C,\mathfrak{i}-}^{k}\}_{k\in\mathcal{V}_{\delta,C,\mathfrak{i}}}$, 
and $\{\eta_{\mathfrak{i}}^{k}\}_{k\in\bar{\mathcal{V}}_{\delta,\mathfrak{i}}}$ as a reordering of 
$\{l_{C,\mathfrak{i}-}^{k}\}_{k\in\bar{\mathcal{V}}_{\delta,C,\mathfrak{i}}}$.  
In all cases, the basis vectors are normalized such that 
$\| \eta_{\mathfrak{i}}^{k} \|_{\infty} = 1$ for every $\mathfrak{i}$ and $k$.

\begin{remark}
In the diffusion regime, the asymptotic analysis of the RTE \cite{jin2009uniformly,wang2022uniform} shows that the smallest-eigenvalue eigenvectors are nearly parallel to the all-ones vector $\mathbf{e}$. Based on the basis selection criterion, these eigenvectors are placed in $\mathcal{V}_{\delta,C_{-},\mathfrak{i}}$ and $\mathcal{V}_{\delta,C_{+},\mathfrak{i}}$, creating a high degree of linear dependence within $\{l_{C_{-}}^{k_{-}},l_{C_{+}}^{k_{+}}\}_{k_{-}\in\mathcal{V}_{\delta,C_{-},\mathfrak{i}},k_{+}\in\mathcal{V}_{\delta,C_{+},\mathfrak{i}}}$. To maintain numerical stability, we therefore use the orthonormal basis $\{\eta_{\mathfrak{i}}^{k}\}_{k\in\mathcal{V}_{\delta,\mathfrak{i}}}$ instead of the original set.
\end{remark}

Finally, we define projection operators in the finite-dimensional space. For any $\mathfrak{i}\in\mathfrak{I}$ and $k\in\mathcal{V}_{\mathfrak{i}}$, the operator $\mathcal{P}_{\delta,\mathfrak{i}}^{k}$ projects a vector onto the coefficient of its $k$-th angular mode component, defined as:
\begin{equation}
\begin{aligned}
\mathcal{P}_{\delta,\mathfrak{i}}^{k}: \quad & \mathbb{R}^{|\mathcal{V}_{\mathfrak{i}}|} \longrightarrow \mathbb{R} \\
& \forall v \in \mathbb{R}^{|\mathcal{V}_{\mathfrak{i}}|}, \ \exists \{\beta_{\mathfrak{i}}^{k}\}_{k\in\mathcal{V}_{\mathfrak{i}}}, \ \text{s.t.} \ v = \sum_{k\in\mathcal{V}_{\mathfrak{i}}}\beta_{\mathfrak{i}}^{k}\eta_{\mathfrak{i}}^{k}, \quad \mathcal{P}_{\delta,\mathfrak{i}}^{k}v = \beta_{\mathfrak{i}}^{k}.
\end{aligned}
\end{equation}

\subsubsection{Continuity and Boundary Conditions}
\label{subsubsec:continuity_and_boundary_conditions}
Based on the notation introduced above, we derive the corresponding conditions for the ATFPS from the continuity and boundary conditions of the original TFPS. These new conditions act only on the coefficients of the adaptive basis functions, thereby achieving dimension reduction of the system.

First, according to the TFPS continuity condition \eqref{eq:2DTFPS_continuity_cond}, for any interior cell interface $\mathfrak{i}\in\mathfrak{I}_{\mathrm{in}}$, we have:
\begin{equation}
\label{eq:deduct_1}
\begin{aligned}
&\sum_{k\in\mathcal{V}_{\delta, C_{-}}}\alpha_{C_{-}}^{k}\chi_{C_{-}}^{k}(\mathfrak{i}_{\mathrm{mid}}) +\chi_{C_{-}}^{s} + \sum_{k\in\bar{\mathcal{V}}_{\delta, C_{-},\mathfrak{i}}}\alpha_{C_{-}}^{k}l_{C_{-}}^{k} + \sum_{k\in\bar{\mathcal{V}}_{\delta, C_{-},\mathfrak{i}^{\parallel}}\cup \bar{\mathcal{V}}_{\delta, C_{-},\mathfrak{i}^{\perp}}}\alpha_{C_{-}}^{k}\chi_{C_{-}}^{k}(\mathfrak{i}_{\mathrm{mid}})\\
=&\sum_{k\in\mathcal{V}_{\delta, C_{+}}}\alpha_{C_{+}}^{k}\chi_{C_{+}}^{k}(\mathfrak{i}_{\mathrm{mid}}) +\chi_{C_{+}}^{s} + \sum_{k\in\bar{\mathcal{V}}_{\delta, C_{+},\mathfrak{i}}}\alpha_{C_{+}}^{k}l_{C_{+}}^{k} + \sum_{k\in\bar{\mathcal{V}}_{\delta, C_{+},\mathfrak{i}^{\parallel}}\cup \bar{\mathcal{V}}_{\delta, C_{+},\mathfrak{i}^{\perp}}}\alpha_{C_{+}}^{k}\chi_{C_{+}}^{k}(\mathfrak{i}_{\mathrm{mid}}).
\end{aligned}
\end{equation}

Let
\[
\tau_{\delta,\mathfrak{i}} = \sum_{k\in\bar{\mathcal{V}}_{\delta, C_{+},\mathfrak{i}^{\parallel}}\cup \bar{\mathcal{V}}_{\delta, C_{+},\mathfrak{i}^{\perp}}}\alpha_{C_{+}}^{k}\chi_{C_{+}}^{k}(\mathfrak{i}_{\mathrm{mid}}) - \sum_{k\in\bar{\mathcal{V}}_{\delta, C_{-},\mathfrak{i}^{\parallel}}\cup \bar{\mathcal{V}}_{\delta, C_{-},\mathfrak{i}^{\perp}}}\alpha_{C_{-}}^{k}\chi_{C_{-}}^{k}(\mathfrak{i}_{\mathrm{mid}}).
\]
Since $\tau_{\delta,\mathfrak{i}} = O(\delta)$, \eqref{eq:deduct_1} simplifies to:
\begin{equation}
\label{eq:deduct_2}
\begin{aligned}
&\sum_{k\in\mathcal{V}_{\delta, C_{-}}}\alpha_{C_{-}}^{k}\chi_{C_{-}}^{k}(\mathfrak{i}_{\mathrm{mid}}) +\chi_{C_{-}}^{s} + \sum_{k\in\bar{\mathcal{V}}_{\delta, C_{-},\mathfrak{i}}}\alpha_{C_{-}}^{k}l_{C_{-}}^{k} \\
\approx&\sum_{k\in\mathcal{V}_{\delta, C_{+}}}\alpha_{C_{+}}^{k}\chi_{C_{+}}^{k}(\mathfrak{i}_{\mathrm{mid}}) +\chi_{C_{+}}^{s} + \sum_{k\in\bar{\mathcal{V}}_{\delta, C_{+},\mathfrak{i}}}\alpha_{C_{+}}^{k}l_{C_{+}}^{k}.
\end{aligned}
\end{equation}

For any selected basis function index $k\in\mathcal{V}_{\delta,\mathfrak{i}}$, applying the projection operator $\mathcal{P}_{\delta,\mathfrak{i}}^{k}$ to both sides of \eqref{eq:deduct_2} yields:
\begin{equation}
\label{eq:deduct_5}
\sum_{k'\in\mathcal{V}_{\delta, C_{-}}}\alpha_{C_{-}}^{k'}\mathcal{P}_{\delta,\mathfrak{i}}^{k}\chi_{C_{-}}^{k'}(\mathfrak{i}_{\mathrm{mid}}) + \mathcal{P}_{\delta,\mathfrak{i}}^{k}\chi_{C_{-}}^{s} \approx \sum_{k'\in\mathcal{V}_{\delta, C_{+}}}\alpha_{C_{+}}^{k'}\mathcal{P}_{\delta,\mathfrak{i}}^{k}\chi_{C_{+}}^{k'}(\mathfrak{i}_{\mathrm{mid}}) + \mathcal{P}_{\delta,\mathfrak{i}}^{k}\chi_{C_{+}}^{s}.
\end{equation}

Strengthening to an equality constraint gives the ATFPS continuity condition:
\begin{equation}
\label{eq:2DAdaptiveTFPS_continuity_cond}
\sum_{k'\in\mathcal{V}_{\delta, C_{-}}}\alpha_{\delta,C_{-}}^{k'}\mathcal{P}_{\delta,\mathfrak{i}}^{k}\chi_{C_{-}}^{k'}(\mathfrak{i}_{\mathrm{mid}}) + \mathcal{P}_{\delta,\mathfrak{i}}^{k}\chi_{C_{-}}^{s} = \sum_{k'\in\mathcal{V}_{\delta, C_{+}}}\alpha_{\delta,C_{+}}^{k'}\mathcal{P}_{\delta,\mathfrak{i}}^{k}\chi_{C_{+}}^{k'}(\mathfrak{i}_{\mathrm{mid}}) + \mathcal{P}_{\delta,\mathfrak{i}}^{k}\chi_{C_{+}}^{s}.
\end{equation}

Similarly, for any boundary interface $\mathfrak{i}=C\cap\partial\Omega\in\mathfrak{I}_{\mathrm{b}}$ and $k\in\mathcal{V}_{\delta,\mathfrak{i}}$, we obtain the new boundary condition:
\begin{equation}
\label{eq:2DAdaptiveTFPS_boundary_cond}
\sum_{k'\in\mathcal{V}_{\delta,C}}\alpha_{\delta,C}^{k'}\mathcal{P}_{\delta,\mathfrak{i}}^{k}\chi_{C,\mathfrak{i}}^{k'}(\mathfrak{i}_{\mathrm{mid}}) + \mathcal{P}_{\delta,\mathfrak{i}}^{k}\chi_{C}^{s} = \mathcal{P}_{\delta,\mathfrak{i}}^{k}\left(\Psi_{\Gamma^{-},\mathfrak{i}-}(\mathfrak{i}_{\mathrm{mid}})\right).
\end{equation}

\subsubsection{Compressed linear system}
\label{subsubsec:compressed_linear_system}

Using the selected basis function coefficients $\{\alpha_{\delta,C}^{k}\}_{k\in\mathcal{V}_{\delta,C},C\in\mathfrak{C}_{\mathrm{2d}}}$ as unknowns with the continuity condition \eqref{eq:2DAdaptiveTFPS_continuity_cond} and boundary condition \eqref{eq:2DAdaptiveTFPS_boundary_cond}, we assemble the compressed linear system:
\begin{equation}
\label{eq:2DAdaptiveTFPS_linear_sys}
\bar{A}_{\delta}\bar{\alpha}_{\delta} = \bar{b}_{\delta},
\end{equation}
where $\bar{\alpha}_{\delta} = (\alpha_{\delta,C}^{k})_{k\in\mathcal{V}_{\delta,C},C\in\mathfrak{C}_{\mathrm{2d}}}$.

Solving \eqref{eq:2DAdaptiveTFPS_linear_sys} and substituting into \eqref{eq:2DAdaptiveTFPS_sol_smooth} yields the smooth component $\bar{\psi}_{\delta}$ of the solution. For quantities away from boundary layers, $\bar{\psi}_{\delta}$ provides sufficient accuracy. To capture layer information, we reconstruct the layer component $\tilde{\psi}_{\delta}$ from $\bar{\psi}_{\delta}$ and combine them to obtain the full solution $\psi_{\delta}$ with global accuracy.

\subsubsection{Layer reconstruction}
\label{subsubsec:layer_reconstruction}
Based on the computed smooth component, we now reconstruct the layer information $\tilde{\psi}_{\delta}$ by deriving the layer basis function coefficients.

From Section~\ref{subsubsec:continuity_and_boundary_conditions}, we have $\alpha_{\delta,C}^{k} \approx \alpha_{C}^{k}$ for $k \in \mathcal{V}_{\delta,C}$. Substituting into the original TFPS continuity condition and neglecting $\tau_{\delta,\mathfrak{i}}$ gives:
\begin{equation}
\label{eq:deduct_7}
\bar{\psi}_{\delta}|_{C_{-}}(\mathfrak{i}_{\mathrm{mid}}) + \sum_{k'\in\bar{\mathcal{V}}_{\delta, C_{-},\mathfrak{i}}}\alpha_{\delta,C_{-}}^{k'}l_{C_{-}}^{k'}  
= \bar{\psi}_{\delta}|_{C_{+}}(\mathfrak{i}_{\mathrm{mid}}) + \sum_{k'\in\bar{\mathcal{V}}_{\delta, C_{+},\mathfrak{i}}}\alpha_{\delta,C_{+}}^{k'}l_{C_{+}}^{k'}.
\end{equation}
For any interior interface $\mathfrak{i}=C_{-}\cap C_{+}\in\mathfrak{I}_{\mathrm{in}}$ and $k\in\bar{\mathcal{V}}_{\delta,\mathfrak{i}}$, applying $\mathcal{P}_{\delta,\mathfrak{i}}^{k}$ and using $\mathcal{P}_{\delta,\mathfrak{i}}^{k}l_{C_{\pm}}^{k} = 1$, $\mathcal{P}_{\delta,\mathfrak{i}}^{k}l_{C_{\pm}}^{k'} = 0$ for $k' \neq k$, we obtain:
\begin{equation}
\label{eq:2DAdaptive_TFPS_layer1}
\begin{aligned}
\alpha_{\delta,C_{-}}^{k}
&= \mathcal{P}_{\delta,\mathfrak{i}}^{k}\left(\bar{\psi}_{\delta}|_{C_{+}}(\mathfrak{i}_{\mathrm{mid}})\right) - \mathcal{P}_{\delta,\mathfrak{i}}^{k}\left(\bar{\psi}_{\delta}|_{C_{-}}(\mathfrak{i}_{\mathrm{mid}})\right),\quad k \in \bar{\mathcal{V}}_{\delta,C_{-},\mathfrak{i}},\\
\alpha_{\delta,C_{+}}^{k}
&= \mathcal{P}_{\delta,\mathfrak{i}}^{k}\left(\bar{\psi}_{\delta}|_{C_{-}}(\mathfrak{i}_{\mathrm{mid}})\right) - \mathcal{P}_{\delta,\mathfrak{i}}^{k}\left(\bar{\psi}_{\delta}|_{C_{+}}(\mathfrak{i}_{\mathrm{mid}})\right),\quad k \in \bar{\mathcal{V}}_{\delta,C_{+},\mathfrak{i}}.
\end{aligned}
\end{equation}

Similarly, for a boundary interface $\mathfrak{i}=C\cap\partial\Omega\in\mathfrak{I}_{\mathrm{b}}$ and $k \in \bar{\mathcal{V}}_{\delta,\mathfrak{i}}$, we obtain:
\begin{equation}
\label{eq:2DAdaptive_TFPS_layer2}
\alpha_{\delta,C}^{k} = \mathcal{P}_{\delta,\mathfrak{i}}^{k}\left(\Psi_{\Gamma^{-},\mathfrak{i}-}(\mathfrak{i}_{\mathrm{mid}})\right) - \mathcal{P}_{\delta,\mathfrak{i}}^{k}\left(\bar{\psi}_{\delta}(\mathfrak{i}_{\mathrm{mid}})\right).
\end{equation}

Substituting \eqref{eq:2DAdaptive_TFPS_layer1} and \eqref{eq:2DAdaptive_TFPS_layer2} into \eqref{eq:2DAdaptiveTFPS_sol_layer} yields $\tilde{\psi}_{\delta}$. Combining with $\bar{\psi}_{\delta}$ gives the full solution $\psi_{\delta}$ with global accuracy.

As demonstrated in \cite{song2025adaptive}, we have the following a posteriori error analysis for the full-order TFPS solution and the compressed ATFPS solution with layer reconstruction:

\begin{theorem}[A posteriori error estimate]\label{theo:posteriori}
Let $\bm{\psi}$ denote the solution of the full-order TFPS, and let $\bm{\psi}_{\delta}$ denote the solution of the ATFPS with layer reconstruction. Moreover, let $\alpha_{\delta}$ be the coefficient vector corresponding to the smooth component of $\bm{\psi}_{\delta}$ and $\mathbf{A}_{\delta}$ be the compressed discrete operator in ATFPS. Then the error between $\bm{\psi}$ and $\bm{\psi}_{\delta}$ satisfies the following upper bound:
\begin{equation}
\label{eq:posteriori}
\begin{aligned}
\|\bm{\psi} - \bm{\psi}_{\delta}\|_{\infty} \leq 96M^{2}\left(24MC_{\gamma,g,M,\delta,\infty}\|\bar{A}_{\delta}^{-1}\|_{\infty} + 2\right)C_{\gamma,g,M,\delta,\infty}\delta \|\alpha_{\delta}\|_{\infty}.
\end{aligned}
\end{equation}
where $C_{\gamma,g,M,\delta,\infty}$ is a constant independent of $g$ and $\delta$, and grows no faster than $M^{1/2}$.
\end{theorem}

\section{MgNet: Multigrid Network}
\label{sec:MgNet}
Based on the discretization schemes for the steady-state RTE in $x$-$y$ geometry, we obtain a fully discretized linear system. Classical iterative methods such as GMRES suffer from rapidly increasing iteration counts as the mesh is refined, while standard multigrid methods fail to achieve $O(1)$ V-cycle convergence due to the hyperbolic nature of the RTE. In contrast, the recursive skeleton factorization proposed in \cite{ho2012fast} and applied to RTE in \cite{fan2019fast,ren2019fast,song2026fast} effectively captures the multilevel decomposition of the RTE solution operator. This motivates the use of MgNet \cite{he2019mgnet}, a neural network architecture that explicitly encodes such recursive factorization, making it well-suited for learning the RTE solution operator.

The remainder of this section is organized as follows. We first introduce the structure of MgNet, followed by the design of a physics-informed loss function that incorporates a high-frequency filtering mechanism to enhance convergence and robustness.

\subsection{Architecture of the Proposed Neural Network}
\label{subsec:network_structure}
\begin{figure}[htbp]
    \centering
    \includegraphics[width=\linewidth]{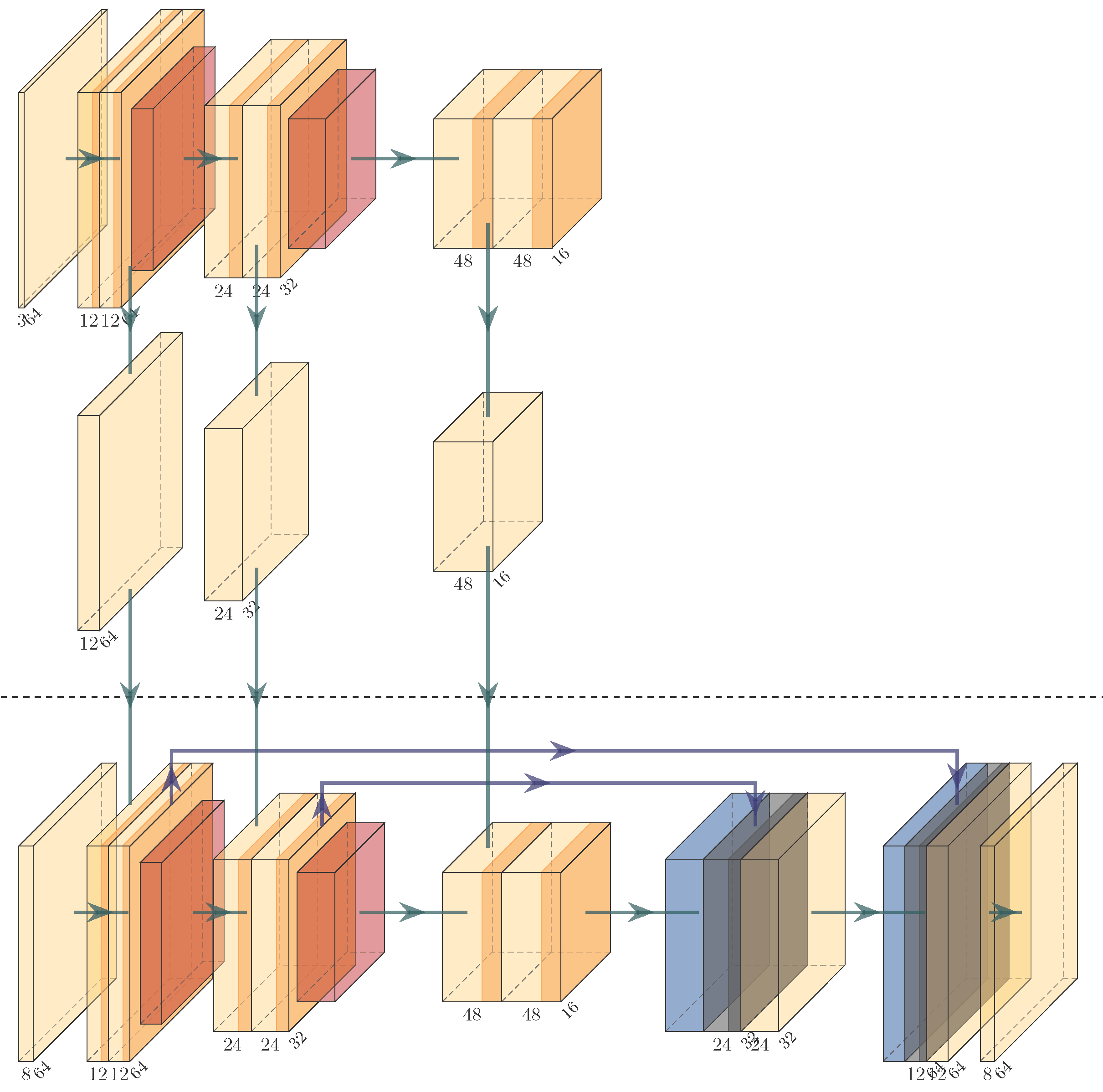}
    \caption{Schematic diagram of CoeffNet and MgNet. The part above the dashed line is CoeffNet, and the part below is MgNet.}
    \label{fig:MixNet}
\end{figure}

The overall network architecture consists of two interconnected sub-networks: \textit{CoeffNet} and \textit{MgNet}. Their collaborative structure enables efficient learning of the solution operator for the steady-state RTE. A schematic diagram of CoeffNet and MgNet is shown in Fig.~\ref{fig:MixNet}, where the part above the dashed line represents CoeffNet and the part below represents MgNet. CoeffNet is coupled with MgNet by passing the medium properties of each layer.

\subsection{CoeffNet: Coefficient Processing Network}
\label{subsec:coeffnet}

The \textit{CoeffNet} sub-model is designed to preprocess and hierarchically distribute the medium properties throughout the computational grid. The network inputs three material coefficients: scaled total cross-section $\sigma_{T}$, scaled absorption cross-section $\sigma_{a}$, and mean free path $\epsilon$. They are concatenated along channels into a 4D tensor $[B, C, H, W]$ with $C=3$ (the three coefficients) and $H,W$ the spatial dimensions. The external source $q$ and boundary conditions are excluded, since the solution operator carries no information about them.

The network processes the input through a series of convolutional layers, batch normalization, and activation functions, organized in a multi-resolution pyramid structure. At each level $l$ of the pyramid, the network produces two outputs:
\begin{itemize}
    \item $\mathbf{a}^{[l]}$: Encoded information about the discretized differential operator at resolution level $l$
    \item $\mathbf{a}^{-1,[l]}$: Encoded information about the corresponding solution operator at resolution level $l$
\end{itemize}

To ensure compatibility with the \textit{MgNet} architecture, the channel dimensions are progressively increased while the spatial dimensions are reduced through strided convolutions. Specifically, for an input tensor of dimensions $[B, C, H, W]$, the outputs at level $l$ have dimensions $[B, C_{l}, H/2^l, W/2^l]$, where $C_{l}=2^{l}C$.
This design is consistent with the structure of the discrete operators and solvers in RSF, and enables the network to learn multi-scale representations of the material properties, which are essential for efficiently solving the RTE across different grid resolutions.

\subsection{MgNet: Multigrid-inspired Neural Network}
\label{subsec:magnet}

The \textit{MgNet} sub-model implements a neural network analog of the classical multigrid V-cycle algorithm. The network takes three inputs:
\begin{enumerate}
    \item The external source term $q$ and the boundary inflow $\Psi_{\Gamma^{-}}$ (discretized on the grid), or equivalently, the right-hand side vector of the fully discretized linear system.
    \item The discretized operator information $\mathbf{a}^{[l]}$ from \textit{CoeffNet} at each level $l$
    \item The solution operator information $\mathbf{a}^{-1,[l]}$ from \textit{CoeffNet} at each level $l$
\end{enumerate}

The network output is the discrete value of angular flux $\psi$ at grid points, which is the solution of the RTE.

The \textit{MgNet} architecture mimics the V-cycle structure of the multigrid method, consisting of pre-smoothing, restriction of residuals from fine to coarse grids, solving on the coarser mesh (performed recursively), prolongation of the coarse-grid correction back to the fine mesh, and post-smoothing. However, it replaces the traditional multigrid components — namely the smoother, restriction operator, and prolongation operator — with the following learnable network structures:

\begin{enumerate}
    \item \textbf{Smoothing operations:} Convolutional layers that update the current solution based on the discrete operator information $a[l]$ and the solution operator information $a^{-1}[l]$ at the current mesh level, analogous to relaxation steps in classical multigrid methods. These layers preserve the number of channels $C$ and the spatial dimensions $H$ and $W$.
    
    \item \textbf{Restriction operators:} Convolutional layers with stride $2$ and kernel size $3$ that transfer residual information from finer to coarser grids. These operations halve the spatial resolution while doubling the channel dimension, enabling the network to represent more abstract, higher-level features.

    \item \textbf{Prolongation operators:} Transposed convolutional layers (or interpolation layers) with stride $2$ and kernel size $4$ that transfer corrected solutions from coarser to finer grids. These operations double the spatial resolution while halving the channel dimension, enabling the network to translate high-level information back to fine-grid details.
    
\end{enumerate}
The tight coupling between coefficient processing and the solution algorithm is a key innovation of our approach. Unlike traditional multigrid methods with fixed operators, our architecture replaces the smoother, restriction, and prolongation operators with learnable neural components, enabling data-driven solution strategies that adapt to material properties. The entire network is trained end-to-end, with \textit{CoeffNet} and \textit{MgNet} jointly optimized to simultaneously learn effective medium representations and tailored solution strategies.

\subsection{loss function structure}
\label{subsec:loss_func}

The proposed loss function is designed in an unsupervised, physics-informed manner, ensuring that the network's output angular flux satisfies the discretized RTE. Specifically, the loss measures the residual of the linear system arising from the discretization scheme.

Let the network's output angular flux be $\psi_{d}$. The loss is defined as the $\ell_2$-norm of the residual:
\begin{equation}\label{eq:loss}
    \mathrm{loss} = \Vert b - A D \psi_{d} \Vert_{2},
\end{equation}
where $D$ maps the discrete angular flux $\psi_{d}$ to the coefficients of the TFPS basis functions. The matrix $A$ and vector $b$ are the system matrix and right-hand side constructed in the standard TFPS scheme ($b$ incorporating source terms and boundary conditions or directly using the network input), as given in Eq. \eqref{eq:2DTFPS_linear_sys}.

To accelerate training and improve robustness, we instead employ the ATFPS as the spatial discretization scheme, which filters out high-frequency modes. This yields the filtered loss:
\begin{equation}\label{eq:loss_delta}
    \mathrm{loss}_{\delta} = \Vert \bar{\mathbf{b}}_{\delta} - \bar{\mathbf{A}}_{\delta} \mathbf{D}_{\delta} \psi_{d} \Vert_{2}.
\end{equation}
Here, $\mathbf{D}_{\delta}$, $\bar{\mathbf{A}}_{\delta}$, and $\bar{\mathbf{b}}_{\delta}$ are defined analogously to their unfiltered counterparts but using the compressed ATFPS basis. Their detailed constructions are provided in Eq. \eqref{eq:2DAdaptiveTFPS_linear_sys}.

\section{Numerical Experiments}
\label{sec:numerical_experiment}
In this section, we conduct comprehensive numerical experiments to validate the proposed model. We evaluate the performance across three distinct physical regimes: diffusion, transport, and sharp interface case. The primary metrics include the convergence behavior of the training and validation losses, and the effectiveness of the trained model as a preconditioner for GMRES solvers.

\subsection{Experimental Configuration}
\label{ssec:experiment_setup}

The physical domain is set to $\Omega = [0,1] \times [0,1]$. The medium parameters $\sigma_T$ and $\sigma_a$ are constructed as piecewise smooth polynomials. Specifically, we utilize a product of one-dimensional polynomials:
\begin{equation}
\Big(2 + \sum_{k=0}^{n} c_{x,k} (x - 0.5)^k \Big) \Big(2 + \sum_{k=0}^{n} c_{y,k} (y - 0.5)^k \Big), \quad n=1,2,3,4,
\end{equation}
where $c_{x,k}, c_{y,k} \sim \mathcal{U}[-1,1]$. The constant term $2$ guarantees the non-negativity of the cross-section. To satisfy the physical constraint $\sigma_a \le \sigma_T$, for each sample we shift $\sigma_T$ such that $\sigma_T(x,y) - \sigma_a(x,y) \ge 0.2$ everywhere in $\Omega$.

The mean free path $\varepsilon$ is varied to simulate different physical regimes:
\begin{itemize}
    \item \textbf{Diffusion regime}: $\epsilon \equiv 0.01(c + 0.1)$, where $c \sim \mathcal{U}[0,1]$.
    \item \textbf{Transport regime}: $\epsilon \equiv 1$.
    \item \textbf{Sharp interface case}: Six configurations mixing diffusion and transport regimes in space, e.g., left/right or top/bottom halves having different $\epsilon$, or a central region ($[0.25,0.75]\times[0.25,0.75]$) being diffusion-dominant while the outer domain is transport-dominant (and vice versa).
\end{itemize}

In our discrete setting, we employ a $32 \times 32$ spatial grid with $4$ discrete velocity directions. These discretization parameters determine the dimensions of the network input (medium parameters) and output (angular flux), thereby defining the specific form of the loss function.
The initial learning rate is set to $0.001$ with a cosine annealing scheduler. For data generation, we produce $10,000$ (or $20,000$) random right-hand sides paired with $1,000$ (or $2,000$) random medium parameters. These are randomly combined to form training and validation samples. The train/validation split is $4:1$, and an independent set of $10$ samples is reserved for testing.

%
%
\subsection{Diffusion Region}
\label{ssec:diffusion_results}

In this example, the dataset consists solely of diffusion-dominated cases ($\varepsilon \ll 1$).

\begin{figure}[htbp]
    \centering
    \includegraphics[width=\linewidth]{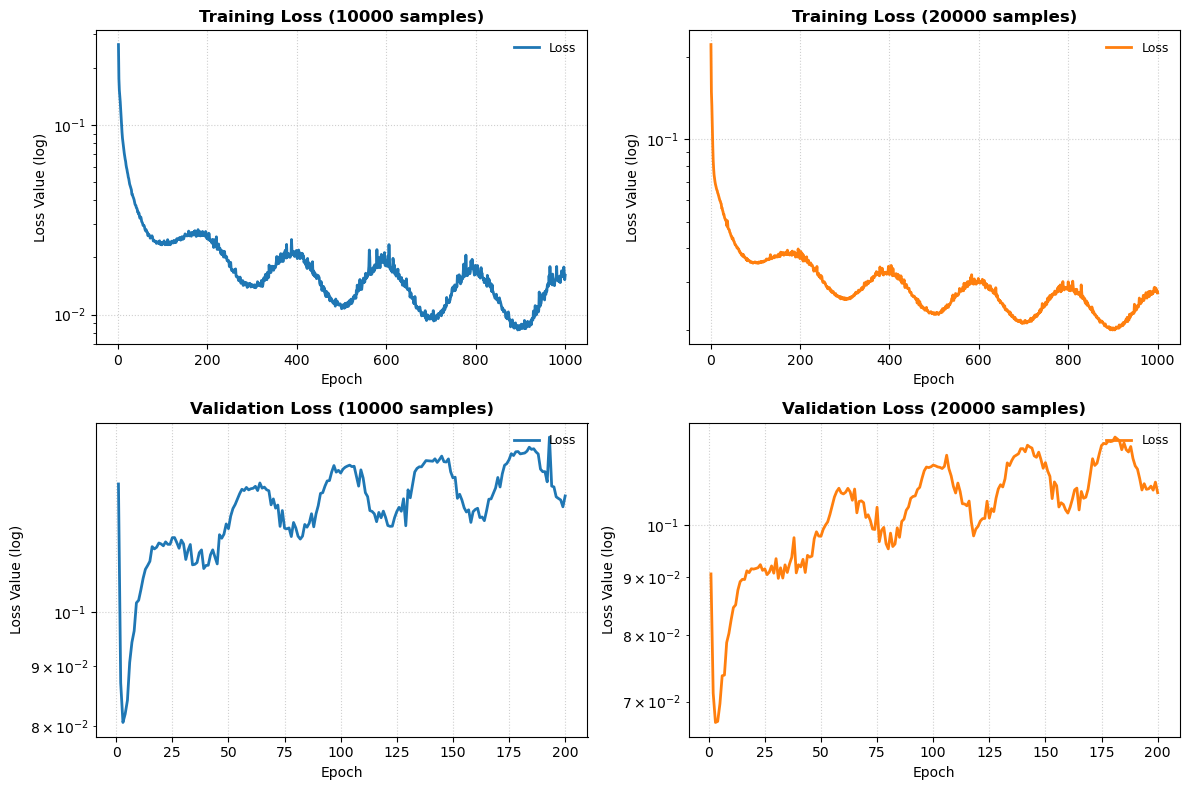}
    \caption{Training and validation loss curves versus epochs for the diffusion region under the full‑order discretization scheme with different sample sizes. }
    \label{fig:diffusion_loss_TFPS}
\end{figure}

Figure \ref{fig:diffusion_loss_TFPS} illustrates the training and validation loss evolution using the full-order TFPS. Although the training loss decreases, the convergence is slow, and the validation loss stagnates or even increases, indicating overfitting. This suggests that the full-order operator contains high-frequency components that are difficult for the network to learn efficiently in the low-diffusion regime.

\begin{figure}[htbp]
    \centering
    \includegraphics[width=\linewidth]{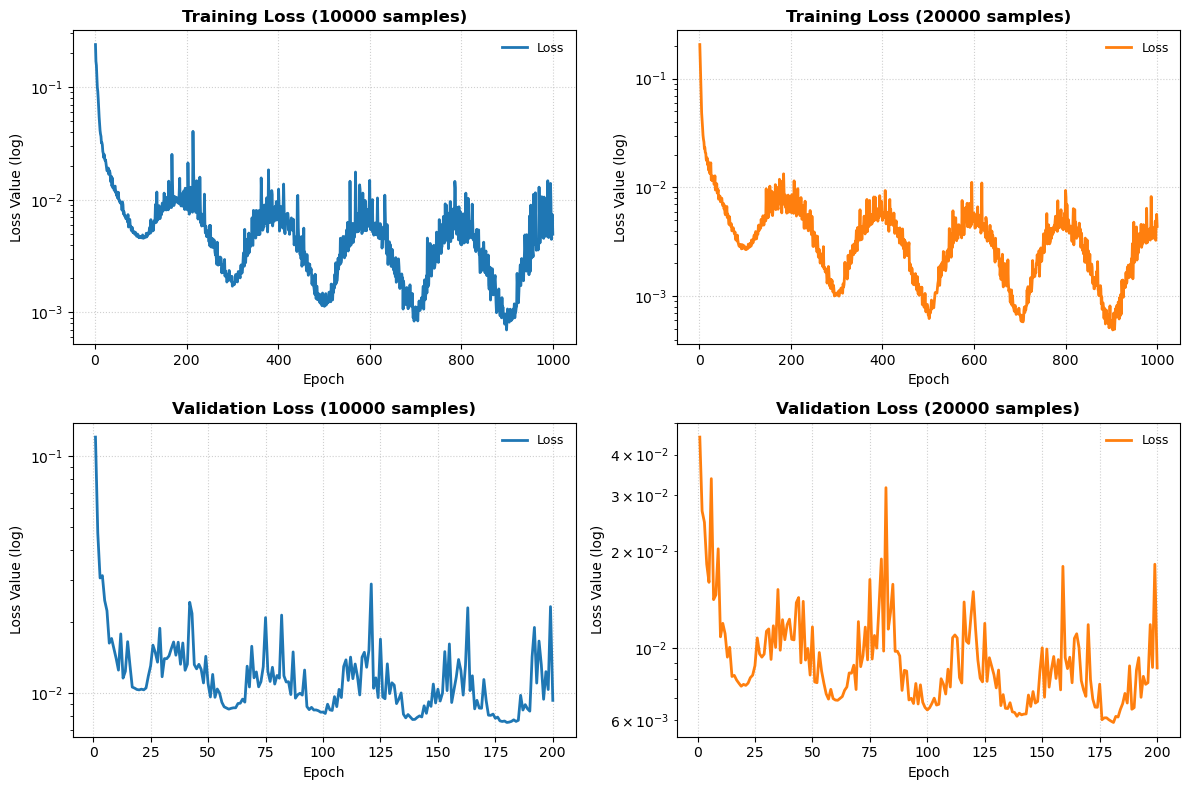}
    \caption{Training and validation loss curves versus epochs for the diffusion region under the compressed discretization scheme with different sample sizes.}
    \label{fig:diffusion_loss_Adaptive_TFPS}
\end{figure}

In contrast, Figure \ref{fig:diffusion_loss_Adaptive_TFPS} presents the results obtained by incorporating the angular domain compression technique into the loss function, i.e., using ATFPS. Both training and validation losses converge smoothly and rapidly. This confirms that by adaptively compressing the operator (discarding rapidly decaying basis functions), the learning difficulty is significantly reduced, allowing the network to capture the underlying low-rank structure of the diffusion operator.

\begin{figure}[htbp]
    \centering
    \includegraphics[width=\linewidth]{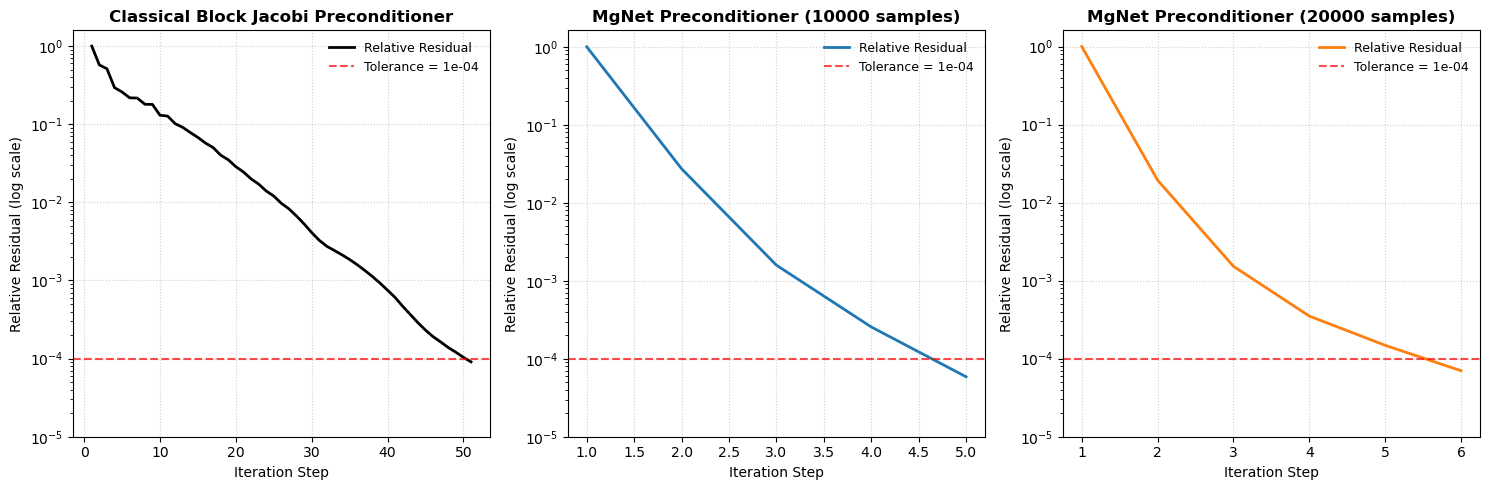}
    \caption{GMRES iteration counts for the test set of the diffusion region problem under different preconditioners.}
    \label{fig:diffusion_test}
\end{figure}

To assess the generalization capability of the trained model, we employ it as a preconditioner for GMRES. Figure \ref{fig:diffusion_test} presents the iteration counts on 10 independent test samples. Compared to the standard Block Jacobi preconditioner, the MgNet-trained preconditioner achieves a dramatic reduction in iteration counts. Interestingly, training with 10,000 samples is sufficient to capture the operator characteristics, as increasing the sample size to 20,000 does not yield further significant improvement.

%
%
\subsection{Transport Region}
\label{ssec:transport_results}

We now consider the pure transport regime where $\epsilon \equiv 1$. In this scenario, the operator lacks the low-rank structure that enables effective adaptive compression. Consequently, ATFPS and TFPS become equivalent, and the two loss functions \eqref{eq:loss} and \eqref{eq:loss_delta} coincide.

\begin{figure}
    \centering
    \includegraphics[width=\linewidth]{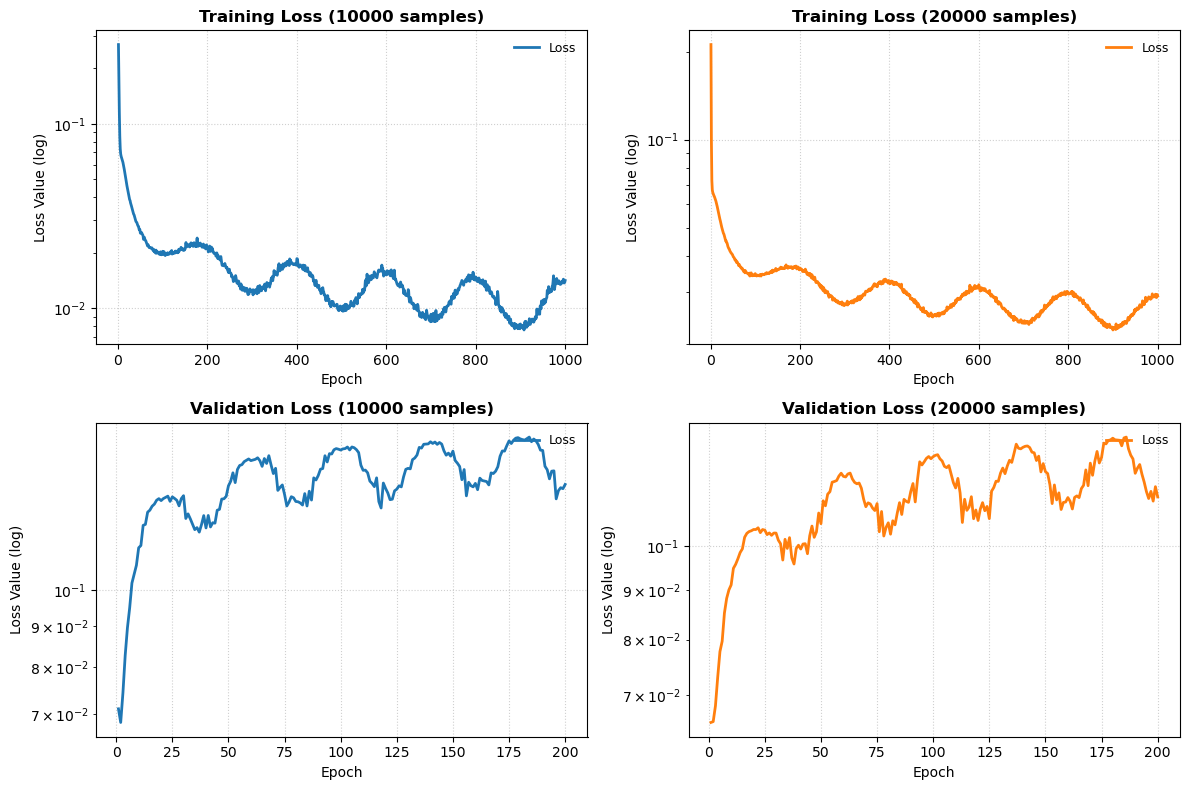}
    \caption{Training and validation loss curves versus epochs for the transport region with different sample sizes. }
    \label{fig:transport_loss}
\end{figure}

Figure \ref{fig:transport_loss} shows the loss curves for this regime. While the training loss exhibits a slight decreasing trend, the convergence is notably slow. Crucially, the validation loss does not follow the training loss; instead, it increases with training epochs. The divergence between training and validation loss clearly indicates a lack of generalization capability. The current network architecture and loss function design are insufficient to capture the complex features of the pure transport operator.

%
\subsection{Sharp Interface}
\label{ssec:interface_results}

This example involves a mixture of transport and diffusion regions, creating sharp interfaces between them (e.g., $\varepsilon=1$ on the left, $\varepsilon=0.01$ on the right). This tests the model's ability to handle multi-scale heterogeneity.

\begin{figure}[htbp]
    \centering
    \includegraphics[width=\linewidth]{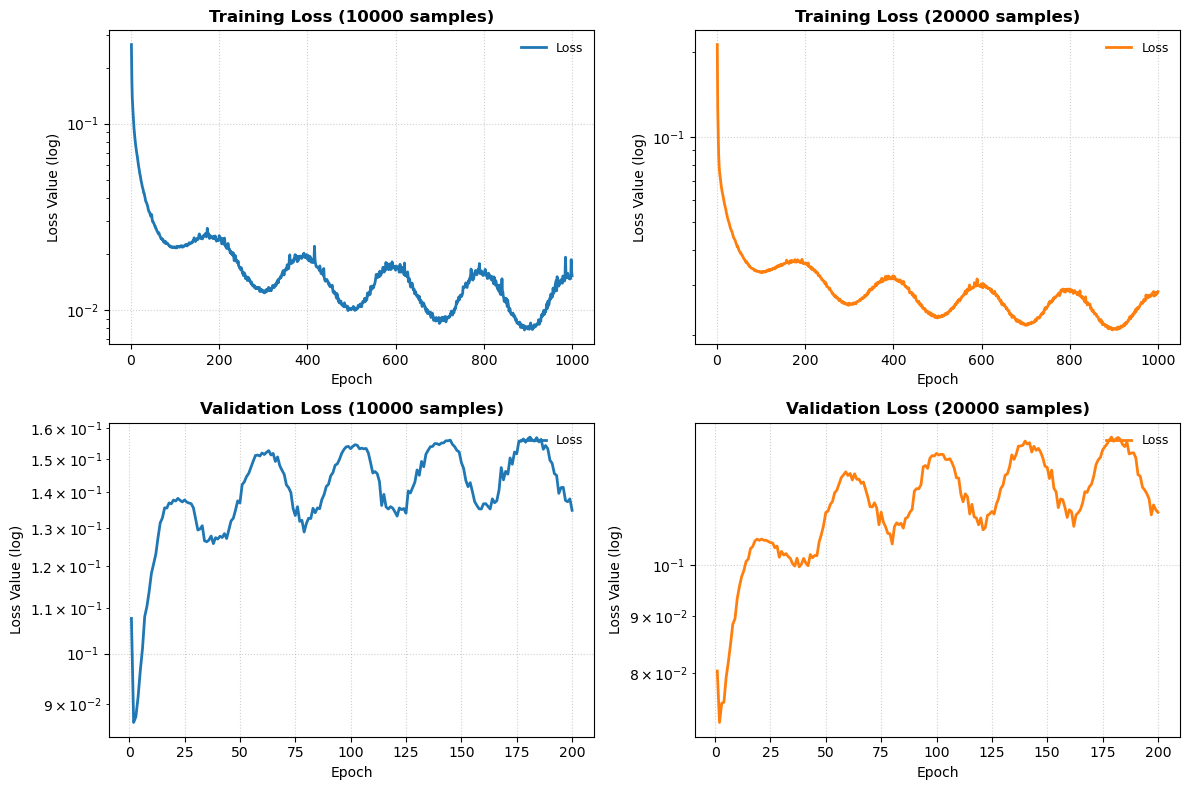}
     \caption{Training and validation loss curves versus epochs for the sharp interface problem under the full‑order discretization scheme with different sample sizes.}
    \label{fig:interface_loss_TFPS}
\end{figure}

\begin{figure}[htbp]
    \centering
    \includegraphics[width=\linewidth]{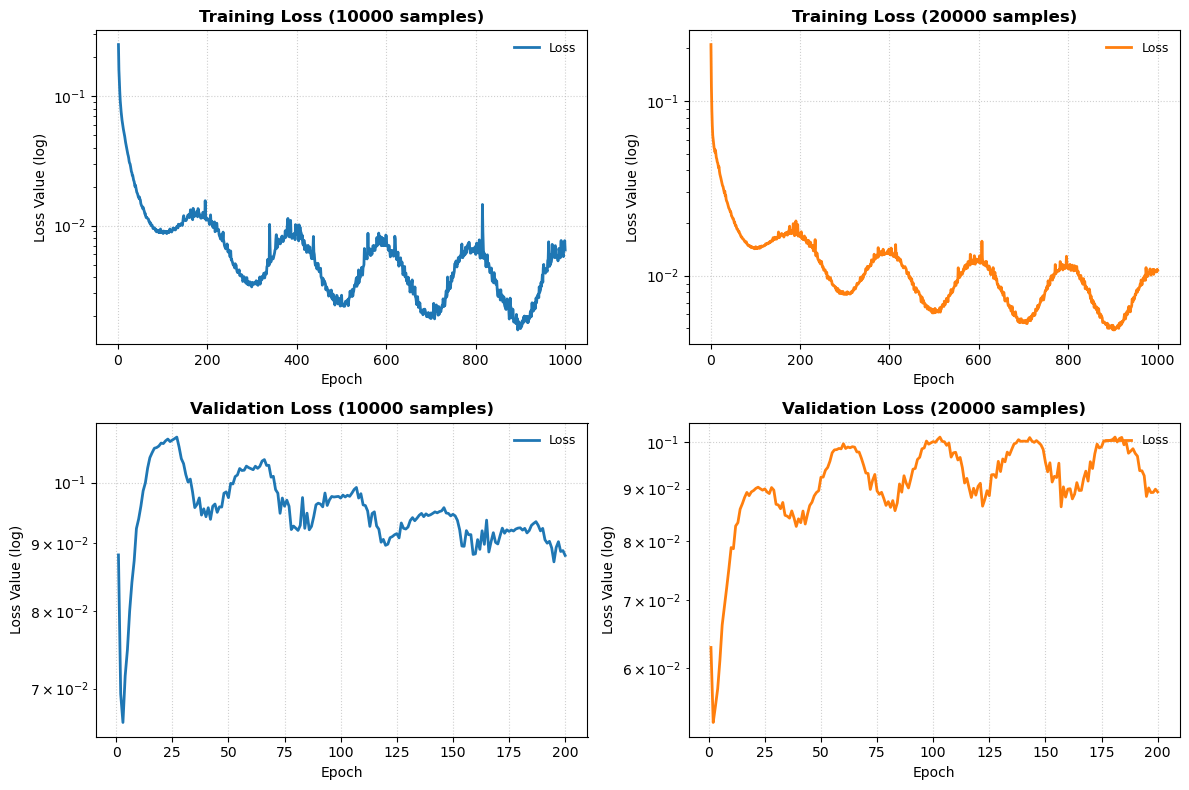}
     \caption{Training and validation loss curves versus epochs for the sharp interface problem under the compressed discretization scheme with different sample sizes.}
    \label{fig:interface_loss_ATFPS}
\end{figure}

The results, shown in Figure \ref{fig:interface_loss_TFPS} and \ref{fig:interface_loss_ATFPS}, are similar to the pure transport case. Despite the presence of diffusion-dominant subregion where learning is easier, the existence of the transport component compromises the overall generalization ability of the model. The validation loss fails to decrease effectively, suggesting that the network struggles to learn a unified representation for the coupled problem with strong discontinuities in $\varepsilon$.

%
\subsection{Conclusion of Numerical Experiments}
\label{ssec:exp_conclusion}

The numerical experiments demonstrate a clear distinction in the performance of the proposed MgNet-based preconditioner:
\begin{enumerate}
    \item \textbf{Diffusion Dominated}: After filtering out the high-frequency modes in the solution space, the model achieves excellent convergence on both training and validation sets. Moreover, when used as a preconditioner, it significantly reduces GMRES iterations compared to traditional methods, demonstrating strong generalization performance.
    \item \textbf{Transport / Sharp Interface}: The model fails to generalize for problems dominated by transport ($\varepsilon \sim 1$) or containing transport regime. The divergence of the validation loss indicates that the current network architecture or loss function is inadequate for capturing the solution operator in such a region.
\end{enumerate}

\section{Summary and Future Work}
\label{sec:summary}
In this work, we adapt the MgNet architecture to solve radiative transfer problems. We first introduce the spatial discretization schemes, including the standard TFPS and the adaptive ATFPS that leverages angular domain compression to filter out high-frequency modes. Then, we incorporate the ATFPS into the loss function design, constructing a physics-informed loss function with high-frequency filtering capability, which is the key innovation of this work. 

The numerical results demonstrate a clear distinction in performance. In diffusion-dominated regimes, the proposed scheme shows excellent convergence and strong generalization when used as a GMRES preconditioner. However, in problems with transport-dominated regimes, the model fails to converge or generalize. Therefore, addressing transport-dominated scenarios remains an important direction for future work, such as designing specialized network architectures or introducing additional regularization terms into the loss function.

\section*{Acknowledgments}
L. Zhang is partially supported by the NSFC grant 12271360, and the Fundamental Research Funds for the Central Universities.

\end{document}